\documentclass[]{article}

\usepackage{mathrsfs, lscape}
\usepackage{amsmath,amsfonts,amssymb,amsthm}
\theoremstyle{plain}
\newtheorem{theorem}{Theorem}[section]
\newtheorem{corollary}[theorem]{Corollary}
\newtheorem{lemma}[theorem]{Lemma}
\newtheorem{proposition}[theorem]{Proposition}
\theoremstyle{definition}
\newtheorem{definition}{Definition}[section]
\theoremstyle{remark}
\newtheorem{remark}{\bf Remark}[section]

\begin{document}
\title{Enumeration of lattices of nullity $k$ and containing $r$ comparable reducible elements}
\maketitle
\begin{center}
\author{Dr.\ A.\ N.\ Bhavale\\hodmaths@moderncollegepune.edu.in\\Head, Department of Mathematics, Modern College of ASC(A),\\ Shivajinagar, Pune 411005, M.S., India.}
\end{center}

\begin{abstract}
In $2002$ Thakare et al.\ counted non-isomorphic lattices on $n$ elements, having nullity up to two. 
In $2020$ Bhavale and Waphare introduced the concept of RC-lattices as the class of all lattices in which all the reducible elements are comparable. 
In this paper, we enumerate all non-isomorphic RC-lattices on $n$ elements. For this purpose, firstly we enumerate all non-isomorphic RC-lattices on $n \geq 4$ elements, having nullity $k \geq 1$, and containing $2 \leq r \leq 2k$ reducible elements. Secondly we enumerate all non-isomorphic RC-lattices on $n \geq 4$ elements, having nullity $k \geq 1$.
This work is in respect of Birkhoff's open problem of enumerating all finite lattices on $n$ elements.
\end{abstract}

MSC Classification $2020$: $06A05,06A06,06A07$.

\section{Introduction}
~\indent
In $1940$ Birkhoff \cite{bir} posed an open problem (which is NP-complete), compute for small $n$ all non-isomorphic posets/lattices on a set of $n$ elements. There were attempts to solve this problem by many authors. In $2002$ Brinkmann and Mckay \cite{bm} counted all non-isomorphic posets with $15$ and $16$ elements. The work of enumeration of all non-isomorphic (unlabeled) posets is still in progress for $n\geq 17$. In the same year, Heitzig and Reinhold \cite{hr} counted non-isomorphic (unlabeled) lattices on up to $18$ elements. 

In $2002$ Thakare et al.\ \cite{tpw} enumerated non-isomorphic lattices on $n$ elements, containing up to $n+1$ edges. These lattices are precisely the lattices having nullity up to two, and contains at most four comparable reducible elements.
Thakare et al.\ \cite{tpw} also enumerated non-isomorphic lattices on $n$ elements, containing two reducible elements.
Recently, Bhavale and Aware \cite{ba} counted non-isomorphic lattices on $n$ elements, containing exactly three reducible elements. 
Observe that, if a lattice contains up to three reducible elements then they are always comparable. Further, if a lattice contains at least four reducible elements then they may not be all comparable. In $2020$ Bhavale and Waphare \cite{bw} introduced the concept of RC-lattices as the class of all non-isomorphic lattices in which all the reducible elements are comparable. Bhavale and Waphare \cite{bw} also introduced the concept of nullity of a poset as the nullity of its cover graph.
In Section \ref{3}, we study the concept of a basic block which is introduced by Bhavale and Waphare \cite{bw}. They have also obtained the number of all non-isomorphic basic blocks of nullity $k$, containing $r$ comparable reducible elements.
In Section \ref{4}, we enumerate all non-isomorphic RC-lattices on $n \geq 4$ elements, using the enumeration of all non-isomorphic basic blocks of nullity $k$, containing $r$ comparable reducible elements. For this purpose, firstly we enumerate RC-lattices on $n \geq 4$ elements, having nullity $k \geq 1$, and containing $r$ reducible elements, where $2 \leq r \leq 2k$. Secondly, we enumerate RC-lattices on $n \geq 4$ elements, having nullity $k \geq 1$.
  
Let $\leq$ be a partial order relation on a non-empty set $P$, and let $(P,\leq)$ be a poset. Elements $x,y \in P$ are said to be $comparable$, if either $x \leq y$ or $y\leq x$. A poset is called a {\it chain} if any two elements in it are comparable. Elements $x,y \in P$ are said to be $incomparable$, denoted by $x||y$, if $x,y$ are not comparable. An element $c\in P$ is {\it{a lower bound}} ({\it{an upper bound}}) of $a,b\in P$ if $c\leq a,c\leq b(a\leq c,b\leq c)$. A {\it{meet}} of $a,b \in P$, denoted by $a \wedge b$, is defined as the greatest lower bound of $a$ and $b$. A {\it{join}} of $a,b\in P$, denoted by $a \vee b$, is defined as the least upper bound of $a$ and $b$. A poset $L$ is a $lattice$ if $a \wedge b$ and $a \vee b$, exist in $L$, $\forall a,b\in L$. Lattices $L_1$ and $L_2$ are {\it{isomorphic}} (in symbol, $L_1\cong L_2$), and the map $\phi:L_1\to L_2$ is an {\it{isomorphism}} if and only if $\phi$ is one-to-one and onto, and $a \leq b$ in $L_1$ if and only if $\phi(a) \leq \phi(b)$ in $L_2$. 
  
An element $b$ in $P$ {\it{covers}} an element $a$ in $P$ if $a<b$, and there is no element $c$ in $P$ such that $a<c<b$. Denote this fact by $a\prec b$, and say that pair $<a, b>$ is a $covering$ or an $edge$. If $a\prec b$ then $a$ is called a {\it{lower cover}} of $b$, and $b$ is called an {\it upper cover} of $a$. An element of a poset $P$ is called {\it{doubly irreducible}} if it has at most one lower cover and at most one upper cover in $P$. Let $Irr(P)$ denote the set of all doubly irreducible elements in the poset $P$. Let $Irr^*(P)= \{x \in Irr(P):x$ has exactly one upper cover and exactly one lower cover in $P\}$. The set of all coverings in $P$ is denoted by $E(P)$. The graph on a poset $P$ with edges as covering relations is called the {\it{cover graph}} and is denoted by $C(P)$. The number of coverings in a chain is called {\it{length}} of the chain.
  
The {\it nullity of a graph} $G$ is given by $m-n+c$, where $m$ is the number of edges in $G$, $n$ is the number of vertices in $G$, and $c$ is the number of connected components of $G$. Bhavale and Waphare \cite{bw} defined {\it nullity of a poset} $P$, denoted by $\eta(P)$, to be the nullity of its cover graph $C(P)$. For $a<b$, the interval $[a,b]=\{x\in P:a\leq x\leq b\}$, and $[a,b)=\{x\in P:a\leq x<b\}$; similarly, $(a,b)$ and ($a,b]$ can also be defined. For integer $n\geq 3$, {\it{crown}} is a partially ordered set $\{x_1,y_1,x_2,y_2,x_3,y_3,\ldots,x_n,y_n\}$ in which $x_i\leq y_i$, $y_i\geq x_{i+1}$, for $i=1,2,\ldots,n-1$, and $x_1\leq y_n$ are the only comparability relations. 
   
An element $x$ in a lattice $L$ is $join$-$reducible$($meet$-$reducible$) in $L$ if there exist $y,z \in L$ both distinct from $x$, such that $y\vee z=x (y\wedge z=x)$. An element $x$ in a lattice $L$ is $reducible$ if it is either join-reducible or meet-reducible. $x$ is $join$-$irreducible$($meet$-$irreducible$) if it is not join-reducible(meet-reducible); $x$ is $doubly$ $irreducible$ if it is both join-irreducible and meet-irreducible. The set of all doubly irreducible elements in $L$ is denoted by $Irr(L)$, and its complement in $L$ is denoted by $Red(L)$. An {\it{ear}} of a loopless connected graph $G$ is a subgraph of G which is a maximal path in which all internal vertices are of degree $2$ in $G$. {\it{Trivial~ear}} is an ear containing no internal vertices. A {\it{non-trivial ear}} is an ear which is not an edge. A vertex of a graph is called {\it pendant} if its degree is one.
For the other definitions, notation and terminology see \cite{dw} and \cite{gg}.

\section{Dismantlable lattices}

In $1974$, Rival \cite{ir} introduced and studied the class of dismantlable lattices. 
\begin{definition}\cite{ir}
\textnormal{A finite lattice of order $n$ is called {\it{dismantlable}} if there exist a chain $L_{1} \subset L_{2} \subset \ldots\subset L_{n}(=L)$ of sublattices of $L$ such that $|L_{i}| = i$, for all $i$}.
\end{definition}

Following result is due to Kelly and Rival \cite{kr}.
\begin{theorem}\cite{kr}\label{crown}
A finite lattice is dismantlable if and only if it contains no crown.
\end{theorem}

By Theorem \ref{crown}, it is clear that an RC-lattice is dismantlable, since it does not contain a crown.

The concept of {\it{adjunct operation of lattices}}, is introduced by Thakare et al.\ \cite{tpw}. 
Suppose $L_1$ and $L_2$ are two disjoint lattices and $(a, b)$ is a pair of elements in $L_1$ such that $a<b$ and $a\not\prec b$. Define the partial order $\leq$ on $L = L_1 \cup L_2$ with respect to the pair $(a,b)$ as follows: $x \leq y$ in L if $x,y \in L_1$ and $x \leq y$ in $L_1$, or $x,y \in L_2$ and $x \leq y$ in $L_2$, or $x \in L_1,$ $ y \in L_2$ and $x \leq a$ in $L_1$, or $x \in L_2,$ $ y \in L_1$ and $b \leq y$ in $L_1$. 
 
It is easy to see that $L$ is a lattice containing $L_1$ and $L_2$ as sublattices. The procedure for obtaining $L$ in this way is called {\it{adjunct operation (or adjunct sum)}} of $L_1$ with $L_2$. We call the pair $(a,b)$ as an {\it{adjunct pair}} and $L$ as an {\it{adjunct}} of $L_1$ with $L_2$ with respect to the adjunct pair $(a,b)$ and write $L = L_1 ]^b_a L_2$. A diagram of $L$ is obtained by placing a diagram of $L_1$ and a diagram of $L_2$ side by side in such a way that the largest element $1$ of $L_2$ is at lower position than $b$ and the least element $0$ of $L_2$ is at the higher position than $a$ and then by adding the coverings $<1,b>$ and $<a,0>$. This clearly gives $|E(L)|=|E(L_1)|+|E(L_2)|+2$. A lattice $L$ is called an {\it{adjunct of lattices}} $L_0, L_1, L_2, \ldots, L_k$, if it is of the form
 $L = L_0 ]^{b_1}_{a_1} L_1 ]^{b_2}_{a_2} L_2 \ldots ]^{b_{k}}_{a_{k}} L_k$. 

Following results are due to Thakare et al.\ \cite{tpw}.
\begin{theorem}\cite{tpw}\label{dac} A finite lattice is dismantlable if and only if it is an adjunct of chains.
\end{theorem}

\begin{corollary}\cite{tpw}\label{k+2}
A dismantlable lattice with $n$ elements has $n+r-2$ edges if and only if it is an adjunct of $r$ chains.
\end{corollary}

Using Corollary \ref{k+2}, Bhavale and Waphare \cite{bw} obtained the following result.
\begin{theorem} \cite{bw} \label{r-1} 
A dismantlable lattice with $n$ elements has nullity $r-1$ if and only if it is an adjunct of $r$ chains.
\end{theorem}

In the following, we prove that a dismantlable lattice has an adjunct representation into chains starting with a maximal chain.
	
\begin{theorem} \label{Cap}
Let $L$ be a dismantlable lattice, and $C$ be a maximal chain in $L$. Then there exist chains $C_1, C_2, \ldots , C_k$ in $L$ such that 
$L = C ]_{a_1}^{b_1} C_1 ]_{a_2}^{b_2} C_2 \cdots ]_{a_k}^{b_k} C_{k}$.
\end{theorem}
\begin{proof} 
We prove the result using induction on $|L| \geq 1$. If $|L|=1$ then we are done. 
Now suppose $|L| > 1$ and the result is true for all dismantlable lattices $K$ with $|K| < |L|$. 
Let $x \in Irr(L)$ and let $L' = L \setminus \{x\}$.
Now we have the following three cases.\\
Case 1: Suppose $x \in C$ and $C  \setminus \{x\}$ is a maximal chain in $L'$.
Now $L'$ is a dismantlable lattice. Therefore by induction hypothesis, we get chains $C_1, C_2, \ldots , C_k$ in $L'$ such that 
$L' = (C  \setminus \{x\}) ]_{a_1}^{b_1} C_1 ]_{a_2}^{b_2} C_2 \cdots ]_{a_k}^{b_k} C_{k}$.
Hence $L = C ]_{a_1}^{b_1} C_1 ]_{a_2}^{b_2} C_2 \cdots ]_{a_k}^{b_k} C_{k}$ as required.\\
Case 2: Suppose $x \in C$ and $C  \setminus \{x\}$ is not a maximal chain in $L'$.
Then there exist $a, b \in C$ such that $a \prec x \prec b$, and a maximal chain $C_0$ in $[a, b]$ such that $x \notin C_0$. But then 
$C' = (C \cap [0, a]) \cup C_0 \cup (C \cap [b, 1])$ is a maximal chain in $L'$. Therefore by induction hypothesis, we get chains $C'_1, C'_2, \ldots , C'_k$ in $L'$ such that 
$L' = C' ]_{a'_1}^{b'_1} C'_1 ]_{a'_2}^{b'_2} C'_2 \cdots ]_{a'_k}^{b'_k} C'_{k}$. 	But then it is easy to see that, 
$L = C ]_a^b C_0 ]_{a'_1}^{b'_1} C'_1 ]_{a'_2}^{b'_2} C'_2 \cdots ]_{a'_k}^{b'_k} C'_{k}$ as required, since $C ]_a^b C_0 = C' ]_a^b \{x\}$.\\
Case 3: Suppose $x \notin C$.
Then $C$ remains a maximal chain in $L'$. Therefore by induction hypothesis, we get chains $C''_1, C''_2, \ldots , C''_k$ in $L'$ such that 
$L' = C ]_{a''_1}^{b''_1} C''_1 ]_{a''_2}^{b''_2} C''_2 \cdots ]_{a''_k}^{b''_k} C''_{k}$.
As $x \notin C$, there exist $a,b\in L$ such that $a\prec x\prec b$. 
If there is an element $z\ne x$ such that $a\prec z\prec b$ then 
$L = C ]_{a''_1}^{b''_1} C''_1 ]_{a''_2}^{b''_2} C''_2 \cdots ]_{a''_k}^{b''_k} C''_{k} \\ ]_a^b \{x\}$.
If there is no $z$ in $L$ such that $z \neq x$ and $a \prec z \prec b$, then we must have $C''_i$ containing both $a$ and $b$. 
Replacing $C''_i$ by $C'''_i = C''_i \cup \{x\}$ we get the required result.
\end{proof}

In the following result, we obtain the bounds on the number of reducible elements of a dismantlable lattice depending on the nullity.
\begin{lemma} \label{r2k}
If $L$ is a lattice of nullity $k \geq 1$ then $2 \leq |Red(L)| \leq 2k$.
\end{lemma}
\begin{proof}
Let $L$ be a lattice of nullity $k$, containing $r$ reducible elements. Therefore by Theorem \ref{r-1}, $L$ is an adjunct of $k+1$ chains. Hence an adjunct representation of $L$ consists of $k$ adjunct pairs, say $(a_i, b_i), \; 1 \leq i \leq k$. Now $Red(L) = \{ a_i, b_i | 1 \leq i \leq k \}$. 
Thus $2 \leq |Red(L)| \leq 2k$, since these reducible elements may not all be distinct.
\end{proof}

\section{Basic blocks} \label{3}
~\indent
Thakare et al.\ \cite{tpw} have defined a \textit{block} as a finite lattice in which the largest element is join-reducible and the least element is meet-reducible.

If $M$ and $N$ are two disjoint posets, the {\it{direct sum}} ({\it{see \cite{rp}}}), denoted by $M \oplus N$, is defined by taking the following order relation on $M\cup N:x\leq y$ if and only if $x,y\in M$ and $x\leq y$ in $M$, or $x,y\in N$ and $x\leq y$ in $N$, or $x\in M,y\in N$. 
In general, $M \oplus N \neq N \oplus M$. Also, if $M$ and $N$ are lattices then $|E(M\oplus N)|=|E(M)|+|E(N)|+1$.

\begin{remark}
Let $L$ be a finite lattice which is not a chain. Then $L$ contains a unique maximal sublattice which is a block, called the {\it{maximal block}}. The lattice $L$ has the form either $C_1\oplus\textbf{B}$ or $\textbf{B}\oplus C_2$ or $C_1\oplus\textbf{B}\oplus C_2$, where $C_1, C_2$ are the chains, and $\textbf{B}$ is the maximal block. Further, $\eta(L) = \eta(\textbf{B})$ and $Red(L) = Red(\textbf{B})$.
\end{remark}

\begin{definition}\label{basicblock}
\textnormal{A block $B$ is a {\it{basic block}} if it is one element or $Irr(B) = \phi$ or removal of any doubly irreducible element reduces nullity by one}.
\end{definition}

In the following result, we obtain the characterization of basic blocks of nullity $k$ in which the reducible elements are all comparable.
\begin{theorem} \label{chbb} 
Let $B$ be an RC-lattice on $n$ elements. Then $B$ is a basic block having nullity $k$ if and only if 
$B = C_0 ]^{b_1}_{a_1} C_1 ]^{b_2}_{a_2} C_2 \cdots ]^{b_{k}}_{a_{k}} C_{k}$ with $a_i, b_i \in C_0$, satisfying 
(i) $|C_i| = 1$, for all $i, \; 1 \leq i \leq k$, (ii) $n = |C_0| + k$,  (iii) $Irr(B) = k+m$, where $m$ is the number of distinct adjunct pairs $(a_i, b_i)$ such that the interval $(a_i, b_i) \subseteq Irr(B)$, and (iv) $|C_0| = |Red(B)| + m$.
Further, if $|Red(B)|=r$ then $n=r+m+k$.
\end{theorem}
\begin{proof}
Suppose $B$ is a basic block having nullity $k$. Therefore by Theorem \ref{r-1}, $B$ is an adjunct of $k+1$ chains. As $B$ is an RC-lattice, all the reducible elements in $B$ are comparable. Therefore by Theorem \ref{Cap}, $B = C_0 ]^{b_1}_{a_1} C_1 ]^{b_2}_{a_2} C_2 \cdots ]^{b_{k}}_{a_{k}} C_{k}$, where $C_0$ is a maximal chain with 
$a_i, b_i \in C_0$, for all $i, \; 1 \leq i \leq k$.
Suppose for some $i, \; 1 \leq i \leq k$, if $|C_i| > 1$ then there exist $x, y \in C_i \cap Irr(B)$, since 
$Red(B) \subseteq C_0$. But then $B \setminus \{x\}$ is a lattice of nullity $k$, a contradiction, since $B$ is a basic block. 
Therefore $|C_i| = 1$ for all $i, \; 1 \leq i \leq k$. Hence $n = |B| = |C_0| + k$. 
Now suppose $(a_{i_1}, b_{i_1}), (a_{i_2}, b_{i_2}), \cdots , (a_{i_m}, b_{i_m})$ are adjunct pairs such that for each 
$j, \; 1 \leq j \leq m$, the interval $(a_{i_j}, b_{i_j}) \subseteq Irr(B)$. Therefore for each $j, \; 1 \leq j \leq m$,
$|(a_{i_j}, b_{i_j}) \cap C_0| = 1$, since $B$ is a basic block.
Also, there is no $x$ such that $x \in C_0 \cap Irr(B)$ but $x \notin (a_i, b_i)$ for all $i, \; 1 \leq i \leq k$.
For if, suppose there is $x \in C_0 \cap Irr(B)$ but $x \notin (a_i, b_i)$ for all $i, \; 1 \leq i \leq k$. 
But then $B \setminus \{x\}$ is of nullity $k$, a contradiction, since $B$ is a basic block. 
Therefore $|Irr(B)| = k + m$, since $|C_i| = 1$ and $C_i \subseteq Irr(B)$ for all $i, \; 1 \leq i \leq k$. 
Also $|B| = |Red(B)| + |Irr(B)|$. But $|B| = |C_0| + k$ implies that $|C_0| + k = |Red(B)| + |Irr(B)|$.
Hence $|C_0| = |Red(B)| + m$, since $|Irr(B)| = k + m$.

Conversely, suppose $B = C_0 ]^{b_1}_{a_1} C_1 ]^{b_2}_{a_2} C_2 \cdots ]^{b_{k}}_{a_{k}} C_{k}$ is a block with $a_i, b_i \in C_0$, satisfying (i) $|C_i| = 1$, for all $i, \; 1 \leq i \leq k$, (ii) $n = |C_0| + k$,  (iii) $Irr(B) = k+m$, where $m$ is the number of distinct adjunct pairs $(a_i, b_i)$ such that the interval $(a_i, b_i) \subseteq Irr(B)$, and (iv) $|C_0| = |Red(B)| + m$.
Now by Theorem \ref{r-1}, nullity of $B$ is $k$. 
Let $C_i = \{ y_i \}$, for all $i, \; 1 \leq i \leq k$. Let $|Red(B)| = r$ and $C_0$ be the chain $x_1 \leq x_2 \leq \cdots \leq x_{r+m}$. Suppose $Red(B) = \{ x_{i_j} | 1 \leq j \leq r \}$. Now if for some $l, \; 1 \leq l \leq k, \; (a_l, b_l) \subseteq Irr(B)$ then 
$|(a_l, b_l) \cap C_0| = 1$, since otherwise $|C_0| - |Red(B)| > m$, which is not possible. Therefore let $(a_l, b_l) \cap C_0 = \{ x_{i_j} \}$, for some $j, \; r+1 \leq j \leq r+m$. But then $a_l = x_{i_j - 1}$ and $b_l = x_{i_j + 1}$. Now $Irr(B) = \{ y_1, y_2, \cdots , y_{k}, x_{i_{r+1}}, \cdots , x_{i_{r+m}} \}$, since $|Irr(B)| = k+m$.
Also for any $z \in Irr(B)$, $B \setminus \{z\}$ is a lattice on $n-1$ elements, having nullity $k-1$, where $n= |C_0| + k$. 
Thus the removal of doubly irreducible element from $B$ decreases its nullity by one. Therefore $B$ is a basic block.
\end{proof}

Bhavale and Waphare \cite{bw} introduced the following concepts namely, retractible element, basic retract, basic block associated to a poset.
\begin{definition}\cite{bw}\label{rtrct} 
\textnormal{Let $P$ be a poset. Let $x \in Irr(P)$. Then $x$ is called a {\it{retractible}} element of $P$ if it satisfies either of the following conditions.}
\begin{enumerate}
\item \textnormal{There are no $y,z \in Red(P)$ such that $y \prec x \prec z$}.
\item \textnormal{There are $y,z \in Red(P)$ such that $y \prec x \prec z$ and there is no other directed path from $y$ to $z$ in $P$.}
\end{enumerate}
\end{definition}

\begin{definition}\cite{bw} \label{brt}
\textnormal{A poset $P$ is a {\it basic retract} if no element of $Irr^{*}(P)$ is retractible in the poset $P$}.
\end{definition}

\begin{definition}\cite{bw}\label{bbas}
A poset $B$ is a {\it basic block associated to a poset} $P$ if $B$ is obtained from the basic retract associated to $P$ by successive removal of all the pendant vertices, and all the retractible elements formed due to removal of all the pendant vertices.
\end{definition}

The following results are due to Bhavale and Waphare \cite{bw}.
\begin{theorem}\cite{bw} \label{redb}
Let $B$ be a basic block associated to a poset $P$. Then $Red(B) = Red(P)$ and $\eta(B) = \eta(P)$.
\end{theorem}

\begin{proposition} \cite{bw} \label{dbb}
If $B_1, B_2$ are basic blocks associated to posets $P_1, P_2$ respectively and $P_1 \cong P_2$ then $B_1 \cong B_2$.
\end{proposition}

\begin{proposition} \cite{bw} \label{ubb}
If $B_1$ and $B_2$ are basic blocks associated to poset $P$ then $B_1 \cong B_2$.
\end{proposition}

In particular, a basic block associated to a lattice is a basic block.

\section{Counting of lattices} \label{4}
~\indent
In this section, we obtain the number of non-isomorphic lattices on $n$ elements in which the reducible elements are all comparable. 
Let ${\mathscr B}(n, k)$ denote the class of all non-isomorphic maximal blocks on $n \geq 4$ elements, having nullity $k \geq 1$, containing reducible elements which are all comparable.
Let ${\mathscr B}(n, k, r) = \{ {\bf B} \in {\mathscr B}(n, k) : |Red({\bf B})| = r \}$.
Let ${\bf B} \in {\mathscr B} (n, k, r)$. Let $B$ be the basic block associated to ${\bf B}$ (see Proposition \ref{ubb}).
Then by Theorem \ref{redb}, $Red(B) = Red({\bf B})$ and $k = \eta(B) = \eta({\bf B})$.
Also $n \geq k+r$ and by Lemma \ref{r2k}, $2 \leq r \leq 2k$.
Let ${\mathscr B}(n; B, k, r) = \{ {\bf B} \in {\mathscr B}(n, k, r) : B$ is the basic block associated to ${\bf B} \}$.
By Theorem \ref{r-1} and Theorem \ref{Cap}, suppose  $B = C_0 ]^{b_1}_{a_1} C_1  ]^{b_2}_{a_2} C_2  \cdots ]^{b_{k}}_{a_{k}} C_{k}$ where $C_0$ is a maximal chain with $a_i, \; b_i \in C_0,$ for all $i, \; 1 \leq i \leq k$. 
By Theorem \ref{chbb}, (i) $|C_i| = 1,$ for all $i, \; 1 \leq i \leq k$, (ii) $n = |C_0| + k$,  (iii) $Irr(B) = k+m$, where $m$ is the number of distinct adjunct pairs $(a_i, b_i)$ such that the interval $(a_i, b_i) \subseteq Irr(B)$, and (iv) $|C_0| = r + m$.
Note that $m \geq 0$. 
Consider a chain $C : x_1 \leq x_2 \leq \cdots \leq x_r$ of the reducible elements of $B$. Clearly $C \subseteq C_0$.
Let $M = \{(x_i, x_{i+1}) \;|\; (x_i, x_{i+1}) = (a_j, b_j)$ for some $j, 1 \leq j \leq k\}$.
Then $m = |M|$ and $|(x_i, x_{i+1}) \cap C_0| = 1$, for all $(x_i, x_{i+1}) \in M$.
As $Red({\bf B}) = Red(B)$ and $|Irr(B)| = m + k$, $|B| = |Red(B)| + |Irr(B)| = r + m + k$.
Therefore $B \in {\mathscr B} (r+m+k, k, r)$ and $n\geq r+m+k$.

Let us denote those $m$ adjunct pairs, if they exist in $B$, by $(a'_j, b'_j)$, $1 \leq j \leq m$. 
Let $m_j$ be the multiplicity of the adjunct pair $(a'_j, b'_j)$ in $B$.
Let $p \geq 0$ be the number of ordered pairs $(x_i, x_j)$, $j > i+1$ such that 
$(x_i, x_j)$ is an adjunct pair in $B$, and $(x_i, x_j) \cap Red(B) \neq \phi$ (that is, the interval $(x_i, x_j)$ contains at least one reducible element).
Let us denote those $p$ adjunct pairs, if they exist in $B$ by $(a''_j, b''_j)$, $1 \leq j \leq p$. 
Let $p_j$ be the multiplicity of the adjunct pair $(a''_j, b''_j)$ in $B$. 
Now we have the following. 
\begin{enumerate}
\item Clearly $k = \displaystyle\sum_{j=1}^{m} m_j + \displaystyle\sum_{j=1}^{p} p_j$ with either $m \geq 1$ or $p \geq 1$, since $k \geq 1$.
\item Hence $m+k = m + \displaystyle\sum_{j=1}^{m} m_j + \displaystyle\sum_{j=1}^{p} p_j = 
\displaystyle\sum_{j=1}^{m} (m_j + 1) + \displaystyle\sum_{j=1}^{p} p_j$.
\item As $|Irr(B)| = m+k$, there are $m+k$ chains of doubly irreducible elements in ${\bf B}$ which correspond to those $m+k$ doubly irreducible elements in $B$.
\item There are $l = r - m - 1$ edges, if they exist, on a maximal chain containing all the reducible elements in $B$. Therefore corresponding to those $l$ edges in $B$, 
there are $l$ edges (or chains of doubly irreducible elements) in ${\bf B}$.
\item As $|{\bf B}| = n$ and $|Red({\bf B})| = r$, $|Irr({\bf B})| = n - r$, and these $n-r$ elements can be spread into $m+p+l$ parts, say $u_i, \; 1 \leq i \leq m+p+l$ satisfying $u_i \geq m_i + 1$ for $1 \leq i \leq m$, $u_{m+i} \geq p_i$ for $1 \leq i \leq p$, and $u_i \geq 0$ for $m+p+1 \leq i \leq m+p+l$.
\item Let $S = \{u = (u_1, u_2, \ldots , u_{m+p+l}) : u_1 + u_2 + \cdots + u_{m+p+l} = n-r$, where $u_i \geq m_i + 1$ for $1 \leq i \leq m$, $u_{m+i} \geq p_i$ for $1 \leq i \leq p$, and $u_i \geq 0$ for $m+p+1 \leq i \leq m+p+l\}$.

\end{enumerate}

Let $P_n^k$ denote the number of partitions of an integer $n$ into $k$ (non-zero, non-decreasing, integer) parts.
Let $B_r(k)$ be the class of all non-isomorphic basic blocks of nullity $k$, containing $r$ comparable reducible elements.
In $2020$ Bhavale and Waphare \cite{bw} enumerated the class $B_r(k)$ (see Theorem $5.13$ of \cite{bw}).
Let $l, m, p, m_i (1\leq i\leq m), p_j (1\leq j\leq p)$ be the integers associated to $B \in B_r(k)$ as defined above, where 
$k = \displaystyle\sum_{j=1}^{m} m_j + \displaystyle\sum_{j=1}^{p} p_j$. By Theorem \ref{chbb}, $|B| = r+m+k$. Then we have the following result.

\begin{theorem} \label{mt}
Let $B \in B_r(k)$, where $k = \displaystyle\sum_{j=1}^{m} m_j + \displaystyle\sum_{j=1}^{p} p_j$ and $2 \leq r \leq 2k$. Let $n \geq |B|=r+m+k$.
Then $$|{\mathscr B}(n; B, k, r)| = \displaystyle\sum_{n-r = u_1 + u_2 + \cdots + u_{m+p+l}} \left( \left(\prod_{i=1}^m P_{u_i}^{m_i + 1} \right) \times \left(\prod_{i=1}^{p} P_{u_{m+i}}^{p_i} \right)\right),$$
where $u_i \geq m_i + 1$ for $1 \leq i \leq m$, $u_{m+i} \geq p_i$ for $1 \leq i \leq p$, and $u_i \geq 0$ for $m+p+1 \leq i \leq m+p+l$.
\end{theorem}
\begin{proof}
Let ${\bf B} \in {\mathscr B}(n; B, k, r)$. Therefore ${\bf B}$ is a maximal block in ${\mathscr B} (n, k, r)$ such that $B$ is the basic block associated to ${\bf B}$. 
By Theorem \ref{redb}, $Red(B) = Red({\bf B})$ and $k = \eta(B) = \eta({\bf B})$. So clearly $B \in B_r(k)$.
As ${\bf B}$ contains $r$ reducible elements, the remaining $n-r$ irreducible elements of ${\bf B}$ can be distributed in $m+p+l$ parts of ${\bf B}$ satisfying the conditions 
$u_i \geq m_i + 1$ for $1 \leq i \leq m$, $u_{m+i} \geq p_i$ for $1 \leq i \leq p$, and $u_i \geq 0$ for $m+p+1 \leq i \leq m+p+l$.
Note that $k = \displaystyle\sum_{i=1}^{m} m_i + \displaystyle\sum_{i=1}^{p} p_i$.
Recall that $S = \{u = (u_1, u_2, \ldots , u_{m+p+l}) : u_1 + u_2 + \cdots + u_{m+p+l} = n-r$, where $u_i \geq m_i + 1$ for $1 \leq i \leq m$, $u_{m+i} \geq p_i$ for $1 \leq i \leq p$, and $u_i \geq 0$ for $m+p+1 \leq i \leq m+p+l\}$.
Then for each $u \in S$, there exists ${\bf B} \in {\mathscr B} (n; B, k, r)$ and vice versa.

Now out of $m+p+l$ parts of $n-r$, the first $m$ parts satisfy $u_i \geq m_i + 1$ for $1 \leq i \leq m$. 
Therefore for fixed $i$, $u_i$ is partitioned into $m_i + 1$ parts, say $u_i^{(1)}, u_i^{(2)}, \ldots , u_i^{(m_i+1)}$ in $P_{u_i}^{m_i + 1}$ ways. 
Further, the next $p$ parts satisfy $u_{m+i} \geq p_i$ for $1 \leq i \leq p$. Therefore for fixed $i$, $u_{m+i}$ is partitioned into $p_i$ parts, say 
$u_{m+i}^{(1)}, u_{m+i}^{(2)}, \ldots , u_{m+i}^{(p_i)}$ in $P_{u_{m+i}}^{p_i}$ ways. 
Furthermore, the remaining $l$ parts satisfying $u_i \geq 0$ for $m+p+1 \leq i \leq m+p+l$ are assigned in a unique way.

Let $S' = \{u' = (u_1^{(1)}, u_1^{(2)}, \ldots , u_1^{(m_1+1)}, u_2^{(1)}, u_2^{(2)}, \ldots , u_2^{(m_2+1)}, \ldots , u_m^{(1)}, u_m^{(2)}, \\
\ldots , u_m^{(m_m+1)}, u_{m+1}^{(1)}, u_{m+1}^{(2)}, \ldots , u_{m+1}^{(p_1)}, u_{m+2}^{(1)}, u_{m+2}^{(2)}, \ldots , u_{m+2}^{(p_2)}, \ldots , u_{m+p}^{(1)}, u_{m+p}^{(2)}, \\
\ldots , u_{m+p}^{(p_p)}, u_{m+p+1}, u_{m+p+2}, \ldots , u_{m+p+l}) : 
0 < u_i^{(1)} \leq u_i^{(2)} \leq \cdots  \leq u_i^{(m_i+1)}, \\
\displaystyle \sum_{j=1}^{m_i+1} u_i^{(j)} = u_i, 1 \leq i \leq m; ~
0 < u_{m+i}^{(1)} \leq u_{m+i}^{(2)} \leq \cdots  \leq u_{m+i}^{(p_i)}, \sum_{j=1}^{p_i} u_{m+i}^{(j)} = u_{m+i}, 1 \leq i \leq p; ~
\sum_{i=1}^{m+p+l} u_i = n-r$, where $u_i \geq m_i + 1$ for $1 \leq i \leq m$, $u_{m+i} \geq p_i$ for $1 \leq i \leq p$, and $u_i \geq 0$ for $m+p+1 \leq i \leq m+p+l\}$.

Then for each $u' \in S'$, there exists ${\bf B} \in {\mathscr B} (n; B, k, r)$ and vice versa.
Further for ${\bf B}, {\bf B'} \in {\mathscr B} (n; B, k, r)$, if there exists $u', v' \in S'$ respectively, then
${\bf B} \cong {\bf B'}$ if and only if $u' = v'$. Therefore ${\mathscr B} (n; B, k, r)$ is equivalent to the set $S'$.
Therefore $|{\mathscr B} (n; B, k, r)| = |S'|$.

Further, $u' \in S'$ if and only if $( \displaystyle \sum_{j=1}^{m_1+1} u_1^{(j)}, \sum_{j=1}^{m_2+1} u_2^{(j)}, \ldots , \sum_{j=1}^{m_m+1} u_m^{(j)}, 
\sum_{j=1}^{p_1} u_{m+1}^{(j)}, \\ \sum_{j=1}^{p_2} u_{m+2}^{(j)}, \ldots , \sum_{j=1}^{p_p} u_{m+p}^{(j)}, u_{m+p+1}, u_{m+p+2}, \ldots , u_{m+p+l} ) \in S$.

Now let $u \in S$. Then for fixed $i, 1 \leq i \leq m$, there are $P_{u_i}^{m_i + 1}$ vectors in $S'$ corresponding to the component $u_i$ of $u$ such that the other components of those vectors of $S'$ other than $u_i^{(1)}, u_i^{(2)}, \ldots , u_i^{(m_i+1)}$ are kept fixed.
Similarly, for fixed $i, 1 \leq i \leq p$, there are $P_{u_{m+i}}^{p_i}$ vectors in $S'$ corresponding to the component $u_{m+i}$ of $u$ such that the other components of those vectors of $S'$ other than $u_{m+i}^{(1)}, u_{m+i}^{(2)}, \ldots , u_{m+i}^{(p_i)}$ are kept fixed.
Thus there are $\displaystyle \left( \prod_{i=1}^m P_{u_i}^{m_i + 1} \right) \times \left( \prod_{i=1}^{p} P_{u_{m+i}}^{p_i} \right)$ vectors in $S'$ corresponding to $u \in S$.
Therefore $|S'| =  \displaystyle\sum_{u \in S} \left( \left( \prod_{i=1}^m P_{u_i}^{m_i + 1} \right) \times \left( \prod_{i=1}^{p} P_{u_{m+i}}^{p_i} \right)\right)$.
Hence $|{\mathscr B} (n; B, k, r)| = |S'| = \displaystyle\sum_{n-r = u_1 + u_2 + \cdots + u_{m+p+l}} \left( \left( \displaystyle\prod_{i=1}^m P_{u_i}^{m_i + 1} \right) \times \left( \displaystyle\prod_{i=1}^{p} P_{u_{m+i}}^{p_i} \right)\right)$,
where $u_i \geq m_i + 1$ for $1 \leq i \leq m$, $u_{m+i} \geq p_i$ for $1 \leq i \leq p$, and $u_i \geq 0$ for $m+p+1 \leq i \leq m+p+l$.
\end{proof}

With the help of the class $B_r(k)$, and using Proposition \ref{dbb} and Theorem \ref{mt}, we obtain in the following result the number of all non-isomorphic maximal blocks on $n$ elements having nullity $k$, containing $r$ reducible elements which are all comparable.
\begin{theorem} \label{mt2}
For $k \geq 1$, for $2 \leq r \leq 2k$, and for $n \geq k+r$,
$$|{\mathscr B}(n, k, r)| = \displaystyle\sum_{B \in B_r(k)} |{\mathscr B}(n; B, k, r)|.$$
\end{theorem}

In the following result using Theorem \ref{mt2} we obtain the number of all non-isomorphic maximal blocks on $n$ elements, having nullity $k$, in which reducible elements are all lying on a chain.
\begin{theorem} \label{mt3}
For $k \geq 1$ and for $n \geq k+3$,
$|{\mathscr B}(n, k)| = \displaystyle\sum_{r=2}^{2k} |{\mathscr B} (n, k, r)|$.
\end{theorem}
\begin{proof}
The proof follows from Lemma \ref{r2k} and the fact that the collection
$\{ {\mathscr B}(n, k, r) : 2 \leq r \leq 2k \}$ forms a partition of ${\mathscr B}(n, k)$.
\end{proof}

In the following result using Theorem \ref{mt3} we obtain the number of all non-isomorphic RC-lattices on $n$ elements, having nullity $k$.
\begin{theorem} \label{mt4}
For $k \geq 1$ and for $n \geq k + 3$, 
$$|{\mathscr L}(n, k)| = \displaystyle\sum_{i=0}^{n - k - 3} (i+1) |{\mathscr B} (n-i, k)|.$$
\end{theorem}
\begin{proof}
It is clear that a lattice $L \in {\mathscr L}(n, k)$ if and only if $L = C \oplus {\bf B} \oplus C'$, where ${\bf B} \in {\mathscr B}(n-i, k)$ and $C, \; C'$ are chains with $|C|+|C'|=i$. For fixed $i \geq 0$, the $i$ elements can be allocated to the chains $C$ and $C'$ in $i+1$ ways. 
Let $j = n - i$. Now for any ${\bf B} \in {\mathscr B}(j, k)$, $j \geq k + 3$. Therefore $i = n - j \leq n - (k + 3) = n - k - 3$. 
Thus $0 \leq i \leq n - k - 3$. Hence the proof.
\end{proof}

Let ${\mathscr L}(n)$ be the class of all non-isomorphic RC-lattices on $n$ elements.
Using Theorem \ref{mt4}, we have the following main theorem.
\begin{theorem} \label{mt}
For $n \geq 1$, $|{\mathscr L}(n)| = 1 + \displaystyle\sum_{k=1}^{n  - 3} |{\mathscr L}(n, k)|$.
\end{theorem}							
\begin{proof}
We know that a chain is the only lattice on $n$ elements of nullity $0$.
Now if $L \in {\mathscr L}(n)$ is a lattice of nullity $k \geq 1$ then $n \geq k + 3$.
Therefore the proof follows from the Theorem \ref {mt4}.
\end{proof}

\noindent
{\bf Conclusion}: In Theorem \ref{mt4} we have obtained the number of all non-isomorphic RC-lattices on $n$ elements, having nullity $k$. Using this enumeration, we have also obtained in Theorem \ref{mt} the number of all non-isomorphic RC-lattices on $n$ elements. We would like to raise the following two problems. \\
{\bf Problem $1$}. Find up to isomorphism all RC-lattices containing $r \geq 4$ reducible elements. \\
{\bf Problem $2$}. Find up to isomorphism all lattices containing $r \geq 4$ reducible elements.


\begin{thebibliography}{10}
\bibitem{ba} A. N. Bhavale and B. P. Aware, {\it Counting of lattices on up to three reducible elements}, The Journal of Indian Math. Soc., (2024), (Accepted).
\bibitem{bw} A. N. Bhavale and B. N. Waphare, {\it{Basic retracts and counting of lattices}}, Australas. J. Combin. \textbf{78(1)}(2020), 73--99.
\bibitem{bir} G. Birkhoff, {\it{Lattice Theory}}, Amer. Math. Soc. Colloq. Pub., 1979.
\bibitem{bm} G. Brinkmann and B. D. McKay, {\it Posets on up to $16$ Points}, Order, \textbf{19}(2002), 147--179.
\bibitem{gg} G. Gr{\"a}tzer. {\it{General Lattice Theory}}, Birkhauser Verlag, 1998.
\bibitem{hr} J. Heitzig and J. Reinhold, {\it{Counting Finite Lattices}}, Algebra universalis, \textbf{48(1)}(2002), 43--53.
\bibitem{kr} D. Kelly and I. Rival, {\it{Crowns, fences and dismantlable lattices}}, Canad. J. Math., \textbf{26(5)}(1974), 1257--1271.
\bibitem{ir} I. Rival, {\it{Lattices with doubly irreducible elements}}, Canad. Math. Bull., \textbf{17}(1974), 91--95.
\bibitem{rp} R. P. Stanley, {\it{Enumerative Combinatorics}}, Wadsworth and Brooks, California, 1986.
\bibitem{tpw} N. K. Thakare, M. M. Pawar, and B. N. Waphare, {\it{A structure theorem for dismantlable lattices and enumeration}}, Periodica Mathematica Hungarica, \textbf{45(1-2)}(2002), 147--160.
\bibitem{dw} D. B. West, {\it{Introduction to Graph Theory}}, Prentice Hall of India, New Delhi, 2002.
\end{thebibliography}
\end{document}